\documentclass[12 pt]{article}
\usepackage{amsmath}
\usepackage{amsthm}
\usepackage{graphicx}

\newcommand{\floor}[1]{\lfloor {#1} \rfloor}

\newcommand{\binomial}[2]{\left( \begin{array}{c} {#1} \\
                        {#2} \end{array} \right)}
\def\la{\langle}
\def\ra{\rangle}

\newtheorem{theorem}{Theorem}[section]
\newtheorem{prop}[theorem]{Proposition}
\newtheorem{lemma}[theorem]{Lemma}

\theoremstyle{remark}
\newtheorem*{remark}{Remark}

\theoremstyle{definition}
\newtheorem{conjecture}{Conjecture}

\title{The Rotor-Router Model}
\author{Lionel Levine \\ University of California, Berkeley}
\date{}

\DeclareSymbolFont{AMSb}{U}{msb}{m}{n}
\DeclareMathSymbol{\N}{\mathbin}{AMSb}{"4E}
\DeclareMathSymbol{\Z}{\mathbin}{AMSb}{"5A}
\DeclareMathSymbol{\R}{\mathbin}{AMSb}{"52}
\DeclareMathSymbol{\C}{\mathbin}{AMSb}{"43}

\begin{document}

\maketitle 

{\bf \large Introduction} \\

Diffusion-limited aggregation (DLA) is a model for dendritic growth in
which a particle executes a random walk in the lattice $\Z^d$ beginning
``at infinity'' and ending when first adjacent to a region $A \subset
\Z^d$.  The terminus of the random walk is then adjoined to the region
$A$, and the procedure is iterated.  It is conjectured \cite{Barlow} that
the resulting region $A$ has fractal dimension exceeding $1$.  Internal
DLA (IDLA) is a variant in which the random walks begin at some fixed
point in $A$ and end when they first leave the $A$.  Unlike DLA, internal
DLA does not give rise to fractal growth.  Lawler et al. (\cite{LBG},
1992) showed that after $n$ walks have been executed, the region $A$,
rescaled by a factor of $n^{1/d}$, approaches a Euclidean ball in
$\R^d$ as $n \rightarrow \infty$.  Lawler (\cite{Lawler95}, 1995)
estimated the rate of convergence.

We propose the following deterministic analogue of IDLA.  First, a cyclic
ordering is specified of the $2d$ cardinal directions in $\Z^d$.  A
``rotor'' pointing in one of the $2d$ cardinal directions is associated to
each point in $A$.  A particle is placed at the origin and is successively
routed in the direction of the rotor at each point it visits until it
leaves the region $A$.  Moreover, every time the particle is routed away
from a given point, the direction of the rotor at that point is
incremented by one step in the cyclic ordering.  When the particle reaches
a point not in $A$, that point is adjoined to the region and given the
first rotor direction in the ordering.

This ``rotor-router'' growth model is similar in many ways to
Engel's ``probabilistic abacus'' (\cite{Engel}, 1975) and to the
abelian sandpile, or ``chip-firing,'' model introduced by Bak et
al{.} (\cite{BTW}, 1988) and studied by Dhar (\cite{Dhar}, 1990)
and Bjorner, et al{.} (\cite{BLS}, 1991).  In the abelian sandpile
model, the points of $A$ are labeled by nonnegative integers,
considered as representing the number of grains of sand at each
point.  A single grain of sand is placed at the origin, and at
each step, every point occupied by at least $2d$ grains of sand
ejects a single grain to each of its $2d$ neighboring lattice
points.  Once the system has equilibrated, with every site
occupied by at most $2d-1$ grains of sand, the process is
repeated. Sandpiles have been studied in part from the point of
view of complex systems, where the interest lies in the
``self-organized criticality'' of the model (cf{.} \cite{BTW,
Creutz, Dhar}).  In contrast to diffusion-limited aggregation,
however, many of the fundamental mathematical properties of which
remain conjectural, sandpiles are tractable.

Much of the mathematics concerning the sandpile model rests on the
abelian nature of the model.  This refers to the nontrivial fact that if
grains of sand are deposited in turn at two different points, allowing the
system to equilibrate both before and after the second grain is deposited,
the resulting configuration does not depend on the order in which the two
grains were deposited.  It is a consequence of a general result of
Diaconis and Fulton (\cite{DF}, 1991) that the rotor-router model has the
same abelian property.  Another result of \cite{DF} implies that
internal DLA is abelian in the sense that the probability that two random
walks with different starting points will terminate at a pair of points
$x$, $y$ is independent of the order in which the walks are performed.

Despite this common abelian property, there are substantive
differences between IDLA and the sandpile model.  While no
asymptotic results are known for the sandpile model, it appears
likely that its asymptotics are not spherical.  In two dimensions,
for example, numerical data indicate that sandpiles may have
polygonal asymptotics.  This is to be contrasted with the main
result of \cite{LBG}, in which IDLA is shown to have spherical
asymptotics.  The rotor-router model may bridge the gap between
sandpiles and IDLA.  Like sandpiles, the rotor-router model is
deterministic; but we conjecture that the asymptotics of the
rotor-router model, like those of IDLA, are spherical.

This paper is intended to serve, first, as a thorough introduction
to the rotor-router model, and second, as a source of conjectures
and open problems which, it is hoped, will inspire future work on
the model. Analogies with IDLA and with the sandpile model are
emphasized throughout.  The paper is structured as follows.

In section 1, we outline some preliminary definitions and prove an
important finiteness lemma.

In section 2, we study the one-dimensional ($d=1$) rotor-router
model in some additional generality.  Fixing positive integers $r$
and $s$, when a particle leaves the interval $A$ through its left
endpoint, $A$ is extended by not just one but $r$ sites to its
left; and when a particle leaves $A$ through its right endpoint,
$A$ is extended by $s$ sites to the right.  We reduce to a set of
states describing the eventual behavior of the system.  By
identifying this collection of states with a suitable subset of
the integer lattice $\Z^3$, we show that the process of adding a
particle at the origin and allowing it to equilibrate can be
realized by a piecewise linear function on $\Z^3$.  We find a set
of invariants of this function that is sufficient to distinguish
between all its orbits.  Our main invariant closely resembles an
invariant of the continuum limit of IDLA in one dimension.

In section 3, we analyze the limiting behavior of the interval $A$
in the generalized one-dimensional model.  Propp (2001)
conjectured that after $t$ iterations, if $A$ is the interval
$[x(t),y(t)]$, then the quantity $\mu(t) = x(t) \sqrt{s} + y(t)
\sqrt{r}$ is bounded independent of $t$. We prove this, and show
furthermore that the limit points of the sequence $\mu$ are
confined to an interval of length $r \sqrt{s} + s \sqrt{r}$, which
is best-possible given that $x$ and $y$ change by increments of
$r$ and $s$, respectively.  We show that in a special case, the
$n$-th particle added ends up on the left or right of the interval
accordingly as the $n$-th term of a certain Sturmian sequence with
quadratic irrational slope is $0$ or $1$.  Finally, we give a
connection with Pythagorean triples.  Given a positive Pythagorean
triple $a^2 + b^2 = c^2$, Propp observed that in the continuum
limit of IDLA, if initially $A$ is the interval $[-a,0]$, then $A$
is the interval $[-c,b]$ at time $a+b-c$, and conjectured that the
same would be true of the rotor-router model in the case $r=s=1$.
We prove this conjecture.

In section 4, we consider the higher-dimensional lattices $\Z^d$,
$d \geq 2$.  We give an estimate for the center of mass of the
region $A$ and prove a weak form of the conjecture that the
limiting shape of the rotor-router model in two dimensions is a
disc.

In section 5, we discuss some outstanding conjectures pertaining to the
one-, two- and three-dimensional models.

\section{\large Preliminaries}

We denote by $\Z$ the set of integers, $\N$ the nonnegative
integers, and $\R$ the real numbers.  If $x \leq y \in \Z$, we
denote by $[x,y]$ the interval $\{z \in \Z : x \leq z \leq y \}$;
we adopt the convention that $[x,y]$ is empty when $x > y$.

Let $e_1, ..., e_d$ be the standard basis vectors for the integer
lattice $\Z^d$, and denote by $E_d = \{\pm e_1, \dots, \pm e_d\}$
the set of cardinal directions in $\Z^d$.  Suppose we are given an
arbitrary total ordering $\leq$ of the set $E_d$: write $E_d =
\{\epsilon_i\}_{i=0}^{2d-1}$ with $\epsilon_0 < \dots <
\epsilon_{2d-1}$.

Let $A \subset \Z^d$ be a finite connected region of lattice
points containing the origin.  A state of the rotor-router
automaton can be described by a pair $(x,l)$, where $x \in A$
represents the position of the particle, and $l : A \rightarrow
[0,2d-1]$ indicates the direction of the rotor at each point.
Define $g(x,l) = (x+\epsilon_{l(x)}, l_x)$, where $l_x$ is the
labeling
    $$ l_x(x') = \begin{cases}
            l(x')+1 \text{ (mod $2n$)} & \text{if } x'=x; \\
            l(x'), & \text{else.}
            \end{cases} $$
This is the state given by routing the particle in the direction of the
rotor $\epsilon_{l(x)}$ and then changing the direction of the rotor.
Write $g^n(0, l) = (x_n, l_n)$.  The sequence of points $(x_0, x_1,
\dots )$ is a lattice path in $A$, possibly self-intersecting, beginning
at the origin.  We denote this path by $p = p(A,l)$.
Lemma~\ref{welldefined}, below, shows that the path $p$ leaves the region
$A$ in finitely many steps.  Let $N$ be minimal such that $x_N \notin A$.
We define $f(A,l) = (A \cup \{x_N\}, l')$, where
    $$ l'(x) = \begin{cases}
            l_N(x) & \text{if } x \in A; \\
            0, & \text{if } x=x_N.
            \end{cases} $$
Then $f$ describes the entire process of adding a particle at the origin
and allowing the system to equilibrate.

\begin{lemma}
\label{welldefined}
The lattice path $p(A,l)$ leaves the region $A$ in finitely many steps.
\end{lemma}

\begin{proof}
If not, the path would visit some point $x \in A$ infinitely many times;
but then it would be routed infinitely many times to each neighbor of $x$.
Inducting along a path from $x$ to a point outside $A$, we conclude the
path does after all leave the region after finitely many steps.
\end{proof}

\section{\large Rotor-router dynamics in one dimension}

In one dimension, rotors alternate between the two directions left
and right; we denote these by $L$ and $R$, respectively.  We
introduce the following generalization of the rotor-router
automaton in one dimension. Let $r$ and $s$ be positive integers.
When the lattice path reaches an unoccupied site $x<0$, all $r$
sites in the interval $[x-r+1, x]$ become occupied; similarly, if
the path reaches an unoccupied site $y>0$, the $s$ sites $[y,
y+s-1]$ become occupied.  In either case, the newly occupied sites
are initially labeled $R$.

Given integers $x \leq 0$ and $y \geq 0$, we denote by ${\bf \Sigma}(x,y)$
the set of all states of the automaton for which the set of occupied sites
is the interval $[x,y]$.  Then ${\bf \Sigma}(x,y)$ is naturally identified
with the set of maps $[x,y] \rightarrow \{R,L\}$.  Let ${\bf \Sigma} =
\bigcup_{x \leq 0 \leq y} {\bf \Sigma}(x,y)$, and denote by $f = f_{r,s} :
{\bf \Sigma} \rightarrow {\bf \Sigma}$ the map on states given by adding
a single particle at the origin and allowing the system to equilibrate.

\subsection{\normalsize Recurrent states}

\begin{lemma}
\label{dumb}
Let $\sigma \in {\bf \Sigma}(x,y)$ be any state.  There exists $N \in \N$
such that $f_{r,s}^N (\sigma) \in {\bf \Sigma}(x',y')$ for $x'<x$, $y'>y$.
\end{lemma}

\begin{proof}
Let $M = 1 + \max(|x|,y)$, and let $N = 2^M$.  After $N$
particles have been deposited and allowed to equilibrate in turn, the
origin will have been visited a total of at least $N$ times.  Since the
rotors at each site alternate pointing left and right, it follows by
induction on $|k|$ that if $|k| \leq M$, then the site $k$ will be visited
at least $2^{M-k}$ times.  Taking $k = \pm M$ proves the lemma.
\end{proof}

Given a state $\sigma \in {\bf \Sigma}(x,y)$, let $0 > u_1 > \dots > u_m$
be the sites to the left of the origin labeled $R$, and let $0 < v_1 <
\dots < v_n$ be the sites to the right of the origin labeled $L$.
Additionally, define $u_0=v_0=0$, $u_{m+1} = x-1$, and $v_{n+1} = y+1$.
Then the path $p(\sigma)$ can be described as follows.

\begin{lemma}
\label{zigzag}
If $\sigma(0) = L$, then the path $p(\sigma)$ travels directly left from
the origin to $u_1$, then right to $v_1$, left to $u_2$, right to $v_2$,
and so on, until it reaches either $u_{m+1}$ or $v_{n+1}$, at which point
it stops.  If $\sigma(0) = R$, the same is true interchanging left with
right and $u_i$ with $v_i$.  In particular, if $m < n$, the path will
terminate at $u_{m+1}$; if $m > n$, it will terminate at $v_{n+1}$; and if
$m = n$, it will terminate at $u_{m+1}$ or $v_{n+1}$ accordingly as
$\sigma(0) = L$ or $\sigma(0) = R$.
\end{lemma}

\begin{proof}
Suppose $\sigma(0) = L$.  Induct on $k$ to show that when the path first
reaches the site $v_k$, it travels left from $v_k$ to $u_{k+1}$, then
right from $u_{k+1}$ to $v_{k+1}$.  When the path first arrives at $v_k$,
it must have previously reached $v_{k-1}$, so by the inductive hypothesis,
the path has come directly right from $u_k$ to $v_k$, hence the entire
interval $[u_k, v_k]$ is now labeled $L$.  Since also, by definition, the
interval $[u_{k+1}+1, u_k - 1]$ is labeled $L$, the path now travels
directly to the left until it reaches the site $u_{k+1}$.  Now the
interval $[u_{k+1}, v_k]$ is entirely labeled $R$, as is the interval
$[v_k+1, v_{k+1}-1]$ by definition; so the path travels directly right
from $u_{k+1}$ to $v_{k+1}$, and the inductive step is complete.

The proof in the case $\sigma(0)=R$ is identical, interchanging the roles
of $u_k$ and $v_k$.
\end{proof}

\begin{remark}
It follows that the $N$ in Lemma~\ref{dumb} can be taken
substantially less than $2^M$.  If, for example, the path
$p\left(f^i (\sigma)\right)$ terminates to the right of the origin
for $i=0,\dots,k$, then for $i \geq 1$ we have $n\left(f^i
(\sigma)\right) = y+(i-1)s$, and so $m\left(f^{i+1}
(\sigma)\right) = m\left(f^i (\sigma)\right) - y - (i-1)s$.  Thus,
to ensure that at least one particle terminates on the left it
suffices to take $N$ large enough so that $Ny + s\frac{N(N-1)}{2}
> |x|$.  Likewise, to ensure that at least one particle terminates
on the right it suffices to have $N|x| + r\frac{N(N+1)}{2} > y$.
Certainly $N = 1+\sqrt{2M/\text{min}(r,s)}$ is enough.
\end{remark}

Let ${\bf Rec}(x,y) \subset {\bf \Sigma}(x,y+s-1)$ be the set of states
for
which there exists an integer $0 \leq i \leq y-x+1$ such that the
interval $[x,x+i-1]$ is entirely labeled $R$, the interval $[x+i,y-1]$ is
entirely labeled $L$, and the interval $[y,y+s-1]$ is entirely labeled
$R$.  Putting $j = y-x-i$, we say as a shorthand that these states are of
the form $R^i L^j R^s$ for nonnegative integers $i$ and $j$.  Let ${\bf
Rec} = \bigcup_{x \leq 0 \leq y} {\bf Rec}(x,y)$.

We are interested primarily in the eventual behavior of $f_{r,s}$ as a
dynamical system on ${\bf \Sigma}$.  The following proposition shows that
for these purposes, it is sufficient to consider the states in ${\bf
Rec}$.  These states will be called the {\it recurrent} states of
$f_{r,s}$.

\begin{prop}
\label{recurrentstates}
The set ${\bf Rec}$ is closed under $f_{r,s}$.  Moreover, for any state
$\sigma \in {\bf \Sigma}$, there is an integer $N \geq 0$ such that
$f_{r,s}^N (\sigma) \in {\bf Rec}$.
\end{prop}

\begin{proof}
Let $\sigma \in {\bf Rec}(x,y)$.  Suppose first that the path $p(\sigma)$
terminates at $y+1$.  Then by Lemma~\ref{zigzag}, $p(\sigma)$ will reach
a site $u<0$ and travel directly right from there
to the unoccupied site $y+1$; in particular, the path never visits the
interval $[x, u-1]$, and $f(\sigma)$ retains its original labels from
$\sigma$ on this interval.  Since $u$ is in the $R^i$ block of $\sigma$,
all of these labels are $R$.  Also, because $p(\sigma)$ travels
to the right from $u$ to $y+1$, the interval $[u, y]$ is entirely labeled
$L$; and the $s$ newly occupied sites in the interval $[y, y+s-1]$ are
labeled $R$.  Thus $f_{r,s}(\sigma)$ has the form $R^i L^j R^s$, as
desired.  On the other hand, if $p(\sigma)$ terminates at $x-1$, then by
Lemma~\ref{zigzag}, it will reach some $v>0$ and travel directly left to
the unoccupied site $x-1$.  In this case, the interval $[x-r, v]$ is
entirely labeled $R$, while the sites in the interval $[v_n+1, y]$ retain
their original labels from $\sigma$.  Since $\sigma \in {\bf Rec}$, the
labels on this interval are of the form $L^j R^s$, so again
$f_{r.s}(\sigma)$ has the desired form.

It remains to show that for any $\sigma \in \Sigma$, some iterate of $f$
takes $\sigma$ into~{\bf Rec}.  By Lemma~\ref{dumb}, there exists $N$ such
that $p\left( f^{N-2}(\sigma) \right)$ terminates on the left and $p\left(
f^{N-1} (\sigma) \right)$ terminates on the right.  Since $p\left(
f^{N-2}(\sigma) \right)$ terminates on the left, we have $f^{N-1}
(\sigma)(t) = R$ for $t \leq 0$.  Since $p\left(f^{N-1} (\sigma)\right)$
terminates on the right, there is some point $u_i < 0$ such that $f^N
(\sigma)$ retains its original labels from $f^{N-1} (\sigma)$ on the
interval $[x,u_i-1]$ and has the form $L^j R^s$ on the interval $[u_i,
y]$; hence $f^N (\sigma)$ is of the form $R^i L^j R^s$ as desired.
\end{proof}

\subsection{\normalsize A piecewise linear function on $\Z^3$}

A state $\sigma \in {\bf Rec}$ has the form $R^i L^j R^s$, and so is
determined by the pair of integers $i$ and $j$ (recall that $s$ was fixed
at the outset, independent of $\sigma$).  If we intend to compute
$f_{r,s}(\sigma)$, however, then the origin must be distinguished so that
we know where to initiate the path $p(\sigma)$.  With the origin
distinguished, ${\bf Rec}$ becomes a three-parameter family of states.
There are several reasonable parameterizations of ${\bf Rec}$, but the one
that simplifies computation most effectively is to let $x,y,z \in \Z$ be
the first occupied site, the first site in the final $R^s$ block, and the
first site in the $L^j$ block, respectively (if $j=0$, we adopt the
convention that $z = y$).  In this way, we identify ${\bf Rec}$ with the
set of integer triples $(x,y,z) \in \Z^3$ satisfying $x \leq 0$, $y \geq
0$ and $x \leq z \leq y$.  Our next proposition determines $f_{r,s}$
explicitly as a piecewise linear function on these triples.

\begin{prop}
\label{piecewise}
$f_{r,s}$ is given on ${\bf Rec}$ by
    $$ f_{r,s} (x,y,z) = \begin{cases}
                   (x,~~~~~~ y+s,~    z-y) & \text{\em if $x+y \leq z$}.
               \\  (x-r,~~    y,~~~~~ z-x+1) & \text{\em if $x+y > z$},

                     \end{cases} $$
\end{prop}

\begin{proof}
Consider first the case $z > 0$, i.e. the case when the origin is
initially labeled $R$.  Then with $m$ and $n$ defined as in
Lemma~\ref{zigzag}, we have $m = -x$ and $n = y-z$.  By
Lemma~\ref{zigzag}, the path $p = p(x,y,z)$ ends up on the right if and
only if $m \geq n$, or $x+y \leq z$.  In this case, once the path reaches
the site $u_n = -n = z-y$, it travels right until it reaches an unoccupied
site.  The first site labeled $L$ is then $z-y$, so
    $$  f(x,y,z) = (x, ~ y+s, ~ z-y).  $$

On the other hand, if $x+y > z$, then the path travels left directly
from the site $v_{m+1} = z+m$ until it reaches an unoccupied site.  If
$v_{m+1} < y-1$, the first site labeled $L$ is then $v_{m+1}+1 = z+m+1 =
z-x+1$; if $v_{m+1} = y-1$, then there are no sites labeled $L$, so
our convention dictates that $z' = y = v_{m+1}+1 = z-x+1$; hence
    $$  f(x,y,z) = (x-r, ~ y, ~ z-x+1). $$

It remains to consider the case $z \leq 0$.  In this case, $m = z-x$ and
$n = y-1$.  The path ends up on the right if and only if $m > n$, or
$x+y \leq z$, as before.  If it ends up on the right, the path travels
directly right from $u_{n+1} = z-n-1 = z-y$ to the unoccupied site, so the
first site labeled $L$ is $z-y$, as desired.  If the path ends up on the
left, then it travels directly left from $v_m = m = z-x$ to an unoccupied
site.  In the case that $v_m < y-1$, the first site labeled $L$ is
then $z-x+1$; if $v_m = y-1$, there are no sites labeled $L$, and by
convention $z' = y = z-x+1$.  This completes the proof.
\end{proof}

We will denote the states $(x, y+s, z-y)$ and $(x-r, y, z-x+1)$ by
$f^+(x,y,z)$ and $f^-(x,y,z)$, respectively.  Thus $f(x,y,z) = f^+(x,y,z)$
or $f^-(x,y,z)$ accordingly as $x+y \leq z$ or $x+y > z$.

\subsection{\normalsize An invariant}

Consider an analogous generalization of stochastic IDLA in one
dimension, in which $r$ or $s$ sites become occupied accordingly
as the random walks terminate on the left or right sight of the
interval.  Each random walk is a ``gambler's ruin'' problem (see,
for example, \cite{GS}), terminating on the right with probability
$\frac{|x|}{|x|+y}$, where $[x,y]$ is the interval of occupied
sites.  Thus the limiting value of the ratio $\frac{|x|}{y}$ as
time goes to infinity satisfies
    $$ \frac{|x|}{y} = \frac{ry/(|x|+y)}{s|x|/(|x|+y)} $$
hence
    $$ \frac{|x|}{y} \rightarrow \sqrt{\frac{r}{s}} ~~~~~~~~~
    \text{as} ~~ t \rightarrow \infty. $$
In fact, in the continuum limit of IDLA, as the frequency with which
particles are dropped is taken to infinity and the interval is rescaled
appropriately, the model becomes deterministic and the quantity
$sx^2-ry^2$ is exactly conserved.  This suggests that the quantity $sx^2 -
ry^2$ is likely to be close to invariant in the rotor-router model as
well.  As the following lemma shows, this is indeed the case.

\begin{lemma}
\label{RRinvariant}
The function
    \begin{equation}
    \label{invariant}
    g_{r,s}(x,y,z) = sx^2 - ry^2 + (r-2)sx + rsy - 2rsz
    \end{equation}
is invariant under $f_{r,s}$.
\end{lemma}

\begin{proof}
Compute
    \begin{eqnarray*}
    g(f^-(x,y,z))-g(x,y,z) &=& g(x-r,y,z-x+1) - g(x,y,z) \\
            &=& s(x-r)^2 - sx^2 - (r-2)sr + 2rs(x-1) \\
            &=& -2rsx + r^2s - r^2s + 2rs + 2rsx - 2rs \\
            &=& 0;
    \end{eqnarray*}
    \begin{eqnarray*}
    g(f^+(x,y,z))-g(x,y,z) &=& g(x,~y+s,~z-y) - g(x,y,z) \\
            &=& ry^2 - r(y+s)^2 + rs^2 + 2rsy. \\
            &=& 0. \qed
    \end{eqnarray*}
\renewcommand{\qedsymbol}{}
\end{proof}

The following section is devoted to showing (Theorem \ref{completeness})
that this invariant, together with the congruence classes of $x$ (mod
$r$) and $y$ (mod $s$), is sufficient to distinguish between all orbits of
$f$.  Later, in sections 3.1-3.3, we use this invariant to derive a
variety of bounds on the growth of the interval $[x,y]$.

\subsection{\normalsize Classification of Orbits}

Since $g$ is linear in $z$, it follows from
Lemma~\ref{RRinvariant} that given any $x,y,n \in \Z$, there is at
most one state $\sigma = (x,y,z) \in {\bf Rec}$ for which
$g(\sigma)=n$; namely, set
    \begin{eqnarray*}
    \label{factoredform}
   z = z_n(x,y) &=& \frac{sx^2-ry^2+(r-2)sx+rsy-n}{2rs} \nonumber \\
        &=& \frac{x(x+r-2)}{2r} - \frac{y(y-s)}{2s} -
        \frac{n}{2rs}
    \end{eqnarray*}
if this is an integer and satisfies $x \leq z \leq y$.

\begin{lemma}
\label{monotonic}
Let $F(x,y)=z_n(x,y)-x$ and $G(x,y)=z_n(x,y)-y$.  Then for all $x\leq
0$, $y \geq 0$ we have
    \begin{eqnarray*}
     (i) &~& F(x,y) \leq F(x-r,y); \\
     (ii) &~& G(x,y) \leq G(x-r,y). \\
     (iii) &~& F(x,y) \geq F(x,y+s);
    \end{eqnarray*}
\end{lemma}

\begin{proof}
From (\ref{factoredform}) we have
    $$ F(x,y) = \frac{x(x-r-2)}{2r} - \frac{y(y-s)}{2s} -
    \frac{n}{2rs}. $$
Thus, for fixed $x$, $F$ is a quadratic in $y$ with maximum at
$y=\frac{s}{2}$; and for fixed $y$, $F$ is a quadratic in $x$ with
minimum at $x=1+\frac{r}{2}$.  This proves (i) and (ii). Likewise,
    $$ G(x,y) = \frac{x(x+r-2)}{2r} - \frac{y(y+s)}{2s} -
    \frac{n}{2rs},$$
so for fixed $y$, $G$ is a quadratic in $x$ with minimum at
$x=1-\frac{r}{2}$, and this proves (iii).
\end{proof}

We adopt the notation $\la x,y,n \ra$ as a shorthand for the state
$(x,y,z_n(x,y))$.  We will say that two states $(x,y,z)$ and
$(x',y',z')$ are {\it congruent} if $x \equiv x'$ (mod $r$) and $y
\equiv y'$ (mod $s$).  By the congruence class of $(x,y,z)$ we
will mean the pair of integer congruence classes $(x ~\text{mod} ~
r, ~~ y ~ \text{mod} ~ s)$.  Trivially, congruence class is
invariant under $f_{r,s}$.   Also, notice that if $z_n(x,y) \in
\Z$ and $\la x',y',n \ra \equiv \la x, y, n \ra$, then $z_n(x',y')
\in \Z$ by (\ref{factoredform}).  Given a state $\sigma_0 = \la
x_0,y_0,n \ra \in {\bf Rec}$ and fixing $y \geq y_0$ congruent to
$y_0$ (mod $s$), let
    $$ A_n(y) = \{x \leq 0 : x \equiv x_0 \text{ (mod $r$)},
            ~ x \leq z_n(x,y) \leq y\}. $$
Then $\la x,y,n \ra \in {\bf Rec}$ if and only if $x \in A_n(y)$.

If $m \in \N$ and $a \leq b \in \Z$ with $a \equiv b$ (mod $m$),
we denote by $[a,b]_m$ the spaced interval $\{a, a+m, a+2m, \dots,
b\}$. This interval is said to have increment $m$.  The following
lemma describes the set of pairs $(x,y)$ for which the triple $\la
x,y,n \ra$ is in ${\bf Rec}$.

\begin{lemma}
\label{spaced} $A_n(y)$ is a nonempty spaced interval with
increment $r$. Furthermore, writing $A_n(y) = [a^-(y), a^+(y)]_r$,
the upper endpoint $a^+$ satisfies $a^+(y+s) \leq a^+(y)$.
\end{lemma}

\begin{proof}
Notice that $x \in A_n(y)$ if and only if
    \begin{equation}
    \label{condition}
    F(x,y) \geq 0 \geq G(x,y).
    \end{equation}
By parts (i) and (ii) of lemma \ref{monotonic}, the set of $x\leq
0$ congruent to $x_0$ (mod $r$) satisfying (\ref{condition}) is a
spaced interval of increment $r$.  This interval is nonempty since
by lemma \ref{dumb} there exists $N \in \N$ such that
$f^N(\sigma_0)$ has right endpoint $y$.

Now since $a^+(y)+r \notin A_n(y)$ and $G(a^+(y)+r,y) \leq
G(a^+(y),y) \leq 0$, it must be that $F(a^+(y)+r,y)<0$. Thus, by
part (iii) of lemma \ref{monotonic}, for any $k \in \N$ we have
    $$ F(a^+(y)+kr,y+s) \leq F(a^+(y)+r, y) < 0, $$
hence $a^+(y+s) \leq a^+(y)$ as desired.
\end{proof}

If ${\bf O} \subset {\bf Rec}$ is an orbit of $f$, we write
$g({\bf O})$ for the constant value of $g$ on {\bf O}, and $c({\bf
O})$ for the common congruence class of the elements of {\bf O}.
Our next result shows that $g$ and $c$ are a complete set of
invariants for $f$ in the sense that no two orbits of $f$ have the
same pair $(g,c)$.

\begin{theorem}
\label{completeness}
Suppose $\sigma, ~ \sigma' \in {\bf Rec}$ are congruent states.  Then
$\sigma$ and $\sigma'$ are in the same orbit of $f_{r,s}$ if and only if
$g_{r,s}(\sigma) = g_{r,s}(\sigma')$.
\end{theorem}

\begin{proof}
The ``only if'' direction is Lemma~\ref{RRinvariant}.  For the
``if'' direction, let $n = g_{r,s}(\sigma)$, and let {\bf O} be
the set of states $(x,y,z) \equiv \sigma$ satisfying $g(x,y,z) =
n$.  If $\alpha=(x,y,z) \in {\bf O}$, then certainly $f(\alpha)
\in {\bf O}$, so {\bf O} contains at least one of the two states
$f^-(\alpha)=(x-r, y, z-x+1)$, $f^+(\alpha) = (x, y+s, z-y)$ in
the definition of the piecewise linear function of
Proposition~\ref{piecewise}.  But $f^+(\alpha)$ is a legal state
if and only if $z-y \geq x$, while $f^-(\alpha)$ is legal if and
only if $z-x+1 \leq y$, and these two conditions are mutually
exclusive.  Thus exactly one of the states $f^+(\alpha)$ and
$f^-(\alpha)$ is in {\bf O}.

By Lemma~\ref{spaced}, for any $y$ we can write
    $$ A_n(y) = [a^-(y), a^+(y)]_r, $$
Denote by $\alpha^{\pm}(y)$ the states $\la a^{\pm}(y), y, n \ra$.
Since $a^-(y)$ is the lower endpoint of the interval $A_n(y)$, we
have $f^-(\alpha^-(y)) \notin {\bf O}$, hence $f^+(\alpha^-(y)) =
\la a^-(y), y+s, n \ra \in {\bf O}$.  Thus
    \begin{equation}
    \label{endpoints}
     a^-(y) \leq a^+(y+s).
    \end{equation}
But by lemma~\ref{spaced}, $a^+(y+s) \leq a^+(y)$, so the state
$\beta := \la a^+(y+s), y, n \ra$ is in {\bf O}. Since
    $$ f^+(\beta) = \alpha^+(y+s) \in {\bf O}, $$
we have $f^-(\beta) \notin {\bf O}$, and it follows that equality
holds in (\ref{endpoints}).  Letting $\sigma_0 \in {\bf O}$ be the
state with minimal $|x|, y$, we conclude that every $\sigma \in
{\bf O}$ is $f^k(\sigma_0)$ for some $k$.
\end{proof}

\section{\large Bounds for the one-dimensional model}

To any state $\sigma \in {\bf Rec}$ we associate an infinite lattice path
$P=P_{r,s}(\sigma)$ in the first quadrant, starting at the origin, whose
$n$-th step is up or to the right accordingly as the $n$-th particle added
to $\sigma$ ends on the right or left of the interval.  This section is
devoted to bounding --- and, in a special case, determining exactly ---
the shape of the path $P$.

\subsection{\normalsize Linear bounds}

The following result shows that the path $P_{r,s}(\sigma)$ is bounded
between two parallel lines of slope $\alpha = \sqrt{\frac{r}{s}}$.  In
particular, for any $\epsilon > 0$, there is a line $\ell$ of slope
$\alpha$ and a translation $\ell'=\ell+(1+\epsilon, -1-\epsilon)$ of
$\ell$ such that all but finitely many steps of $P$ lie between $\ell$ and
$\ell'$ (figure 1).

\begin{figure}
\centering
\includegraphics[scale=.5]{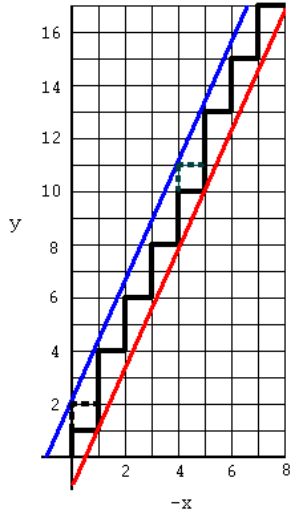}
\caption{\footnotesize The lattice path $P_{5,1}(0,0,0)$ is asymptotically
bounded between two lines of slope $\sqrt{5}$ separated by a translation
of $(1,-1)$.  It crosses the lower line at the points $(1,1), (5,10),
(221,493), (1513,3382), \dots$.  The dotted steps show the unique path
bounded between these two lines, where it differs from $P$.}
\end{figure}

\begin{theorem}
\label{boundedness}
Let $\sigma_0 \in {\bf Rec}$ be any state, and set $\sigma_t =
(x(t),y(t),z(t))$ $= f_{r,s}(\sigma(t-1))$.  Then $|x\sqrt s + y\sqrt r|$
is bounded independent of $t$.  Specifically, given any $\epsilon > 0$
there exists $N \in \N$ such that for all $t > N$,
    \begin{equation}
    \label{bounds}
    -\frac{(r-2)\sqrt{s} - s \sqrt{r} - \epsilon}{2}
        < x \sqrt{r} + y \sqrt{s} <
     \frac{(r+2)\sqrt{s} + s \sqrt{r} + \epsilon}{2}.
    \end{equation}
\end{theorem}

\begin{remark}
The difference between the upper and lower bounds in (\ref{bounds})
approaches $r\sqrt s + s\sqrt r$ as $\epsilon \rightarrow 0$.  This is
best-possible, since the difference of the bounds must be at least the
difference of $x \sqrt{s} + (y+s) \sqrt{r}$ and $(x-r) \sqrt{s} + y
\sqrt{r}$.
\end{remark}

\begin{proof}
Write $u = \frac{r}{2}$, $v = \frac{s}{2}$.  Let $C = g_{r,s}(\sigma)$.
Since $g_{r,s}$ is invariant, we have
    \begin{equation}
    \label{matt}
    sx^2 + (r-2)sx  =  ry^2 - rsy + 2rsz + C
    \end{equation}
at any time $t$.  Then completing the square and using the fact that $z
\leq y$, we obtain
    \begin{eqnarray}
    \label{evan}
    s(x+u-1)^2 - s(u-1)^2  &=& ry^2 - rsy + 2rsz + C \nonumber \\
                   &\leq&  ry^2 + rsy + C  \nonumber \\
                               &=&  r(y+v)^2 - rv^2 + C.
    \end{eqnarray}
Let $C' = s(u-1)^2 - rv^2 + C$.  Lemma~\ref{dumb} ensures that $x
\rightarrow -\infty$ and $y \rightarrow \infty$ as $t \rightarrow \infty$.
In particular, we can take $t$ sufficiently large so that $x \leq 1-u$, $y
\geq v$, and $r(y+v)^2 \geq C'$.  Then we obtain from (\ref{evan})
    \begin{eqnarray}
    \label{adam}
    \sqrt{s} (1-x-u) &\leq& \sqrt{r(y+v)^2 + C'} \\
             &\leq& \sqrt{r}(y+v) + \sqrt{|C'|}, \nonumber
    \end{eqnarray}
which gives a time-independent lower bound for $x\sqrt s + y\sqrt r$.

Similarly, isolating the terms of (\ref{matt}) involving $y$, completing
the square, and using the fact that $z \geq x$ gives
    \begin{eqnarray}
    \label{beth}
    r(y-v)^2 - rv^2 &=& sx^2 + (r-2)sx - 2rsz - C  \nonumber \\
                  &\leq& sx^2 - (r+2)sx - C  \nonumber \\
                     &=& s(x-u-1)^2 - s(u+1)^2 - C.
    \end{eqnarray}
Write $C'' = rv^2 - s(u+1)^2 - C$.  By Lemma~\ref{dumb}, we can take $t$
sufficiently large so that $x \leq -u-1$, $y \geq v$ and $s(x-u-1)^2 \geq
C''$, and estimate
    \begin{eqnarray}
    \label{mira}
    \sqrt{r} (y-v) &\leq& \sqrt{s(x-u-1)^2 + C''} \\
               &\leq& \sqrt{s} (u+1-x) + \sqrt{|C''|} \nonumber,
    \end{eqnarray}
which gives a time-independent upper bound.

To show (\ref{bounds}), choose $t$ sufficiently large so that
$\frac{|C'|}{\sqrt{r}(y+v)}$ and $\frac{|C''|}{\sqrt{s}(u+1-x)}$ are
strictly less than $\epsilon$.  Then using first-order Taylor estimates
for the square roots in (\ref{adam}) and (\ref{mira}), we obtain
  $$ \sqrt{s}(1-x-u) \leq \sqrt{r}(y+v) + \frac{|C'|}{2\sqrt{r}(y+v)}
                   <  \sqrt{r}(y+v) + \frac{\epsilon}{2}; $$
  $$ \sqrt{r}(y-v) \leq \sqrt{s}(u+1-x) + \frac{|C''|}{2\sqrt{s}(u+1-x)}
             <  \sqrt{s}(u+1-x) + \frac{\epsilon}{2}; $$
and hence
    \begin{eqnarray*}
    -\frac{(r-2)\sqrt{s} + s\sqrt{r} + \epsilon}{2}
        &=& (1-u)\sqrt{s} - v\sqrt{r} - \frac{\epsilon}{2} \\
            &<& x \sqrt{s} + y \sqrt{r} \\
        &<& (1+u)\sqrt{s} + v\sqrt{r} + \frac{\epsilon}{2} \\
        &=& \frac{(r+2)\sqrt{s} + s\sqrt{r} + \epsilon}{2}.\qed \\
    \end{eqnarray*}
\renewcommand{\qedsymbol}{}
\end{proof}

\subsection{\normalsize Sturmian words}

For certain values of $r$ and $s$, the inequality (\ref{bounds}) holds for
all $t,t' \in \N$ even when $\epsilon = 0$, and the three coordinates
$x(t), y(t), z(t)$ can be determined exactly in closed form.  In
Proposition~\ref{sturmian}, we treat the case $r=2$, $s=1$.

By a {\it binary word} $w = w_0 w_1 w_2 \dots$ we will mean a map $\N
\rightarrow \{0,1\}$; we write $w_i$ for the image of $i \in \N$ under
this map.  To each state $\sigma \in {\bf Rec}$ we associate a binary word
$w = w_{r,s}(\sigma)$, whose $n$-th term is $0$ or $1$ accordingly as
$f_{r,s}^n(\sigma)$ is $f^-(f^{n-1}\sigma)$ or $f^+(f^{n-1}\sigma)$.  A
word $w = w_0 w_1 w_2 \dots$ is called {\it Sturmian} if it has the form
        \begin{equation}
        \label{floor}
    w_n = \floor{(n+1)\alpha + \beta} - \floor{n\alpha + \beta}.
        \end{equation}
for real numbers $0 \leq \alpha, \beta < 1$, $\alpha$ irrational.
Sturmian words have been extensively studied and have many equivalent
characterizations; see \cite{BS} for a survey.  The word defined by
(\ref{floor}) is called the Sturmian word of slope $\alpha$ and intercept
$\beta$.

\begin{prop}
\label{sturmian}
Suppose $r=2$, $s=1$, and $\sigma=(0,0,0)$.  Then $w_{r,s}(\sigma)$ is
Sturmian with slope $\alpha = \sqrt{2}-1$ and intercept $\beta =
\frac{\alpha}{2}$.
\end{prop}

\begin{proof}
We have $g_{2,1}(\sigma)=0$ and from (\ref{factoredform}),
        \begin{equation}
    \label{innocentlooking}
        z_0(x,y) = \frac{x^2-2y(y-1)}{4}.
        \end{equation}
\noindent We will show that $f^n(\sigma) = \sigma_n :=(x(n),y(n),z(n))$,
where
        $$ x(n) = -2 \floor{(n+\frac12)\alpha}, ~~
           y(n) = n + \frac{x(n)}{2}, ~~
           z(n) = z_0(x(n),y(n)). $$
Since trivially $\sigma_n \equiv \sigma$ and $g_{2,1}(\sigma_n) =
g_{2,1}(\sigma) = 0$, by Theorem~\ref{completeness}, it is sufficient to
show $x(n) \leq z(n) \leq y(n)$.  Note the inequalities
    \begin{equation}
    \label{xobvious}
        -2\alpha(n+\frac12) \leq  x(n) < -2\alpha(n+\frac12) + 2;
    \end{equation}
    \begin{equation}
    \label{yobvious}
        (1-\alpha)n - \frac{\alpha}{2} \leq  y(n)  <
                (1-\alpha)n - \frac{\alpha}{2} + 1.
    \end{equation}
\noindent Since $x^2$ is decreasing on the interval $x<0$, we have from
(\ref{xobvious})
    $$ \alpha^2 \left( n+\frac12 \right)^2
        - 2\alpha \left( n+\frac12 \right) + 1 < \frac{x^2}{4}
    \leq \alpha^2 \left( n+\frac12 \right)^2; $$
and since $y(y-1)$ is increasing on the interval $y>0$, we obtain from
(\ref{yobvious})
    \begin{eqnarray*}
    (1-\alpha)^2 n^2 - (1-\alpha^2)n
        + \frac{\alpha}{2} + \frac{\alpha^2}{4}
    &=& \left( (1-\alpha)n - \frac{\alpha}{2} \right)
        \left( (1-\alpha)n - \frac{\alpha}{2} - 1 \right) \\
    &\leq& y(y-1) \\
    &<& \left( (1-\alpha)n - \frac{\alpha}{2} \right)
        \left( (1-\alpha)n - \frac{\alpha}{2} + 1 \right) \\
    &=& (1-\alpha)^2 n^2 + (1-\alpha)^2 n
        - \frac{\alpha}{2} + \frac{\alpha^2}{4}.
    \end{eqnarray*}
Also notice that $\alpha^2 = 1-2\alpha = \frac12 (1-\alpha)^2$.  Now
(\ref{innocentlooking}) is bounded above by
    \begin{eqnarray*}
    z_0(x,y) &\leq& \alpha^2 \left( n^2 + n + \frac14 \right)
        - \frac{(1-\alpha)^2}{2}n^2 + \frac{1-\alpha^2}{2}n
        - \frac{\alpha}{4} - \frac{\alpha^2}{8} \\
         &=& (1-\alpha)n - \frac{\alpha}{2} + \frac18
        ~~\leq~~ y - \frac{1}{8},
    \end{eqnarray*}
so $z \leq y$ as desired.  Similarly, (\ref{innocentlooking}) is bounded
below by
    \begin{eqnarray*}
    z_0(x,y) &>& \alpha^2 \left( n^2 + n + \frac14 \right)
            -2 \alpha \left( n + \frac12 \right) + 1 \\
~~~~~~~~~~&~&~~~~~~~~~~~~~~
-\frac{(1-\alpha)^2}{2}n^2 - \frac{(1-\alpha)^2}{2}n
            + \frac{\alpha}{4} - \frac{\alpha^2}{8} \\
         &=& -2\alpha n - \alpha + \frac{9}{8}
         ~>~ x - \frac{7}{8}.
    \end{eqnarray*}
Since $x$ and $z_0(x,y)$ are integers, it follows that $z \geq x$.
\end{proof}

For general $r$ and $s$, write $w_{r,s} = w_{r,s}(0,0,0)$.  It is
not true that $w_{r,s}$ is Sturmian for every pair $r,s$.  It is a
classical theorem of Morse and Hedlund (see, e.g., \cite{BS},
\cite{De Luca}) that a Sturmian word has exactly $n+1$ distinct
factors (subwords) of length $n$; and it turns out, for example,
that $w(5,1)$ has $70$ factors of length $68$.  It does not even
appear true that every $w_{r,s}$ is eventually Sturmian.  It does,
however, appear that $w_{r,s}$ is Sturmian for a substantial
number of pairs $(r,s)$.  The set of such pairs is quite complex;
see figure \ref{sturmianregion} and the discussion in section 5.

\subsection{\normalsize Nonlinear bounds and Pythagorean triples}

We now turn to the case $r=s=1$ and consider an initial state consisting
of an interval of occupied sites to the left of the origin.  Our next two
propositions can be seen as rotor-router analogues of the conservation of
$sx^2-ry^2$ in the continuum limit of IDLA (see section 2.3).

\begin{prop}
\label{hyperbolas}
Let $n$ be a positive integer and $\sigma_0 = (-n,0,0)$.  Let
$\sigma_t = (x(t),y(t),z(t)) = f_{1,1}(\sigma_{t-1})$.  Then $x^2 - y^2 <
n^2 + n$ and $(x-1)^2 - (y-1)^2 > n^2$ for all $t \geq 0$.
\end{prop}

\begin{figure}
\centering
\includegraphics[scale=.6]{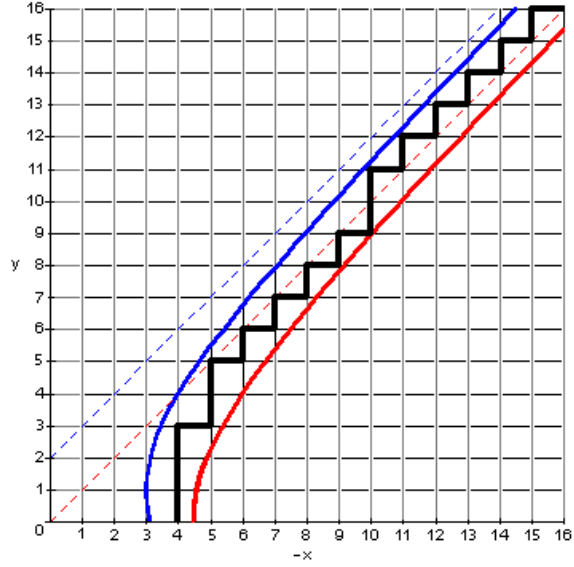}
\caption{\footnotesize The lattice path $P_{1,1}(-4,0,0)$ is bounded
between the hyperbolas $x^2 - y^2 = 20$ and $(x-1)^2 - (y-1)^2 = 16$.}
\end{figure}

\begin{proof}
For the first inequality, we may assume $x^2>y^2$, i.e. $-x>y$.  Then
$x+2z\leq x+2y < y$.  Now from lemma \ref{RRinvariant},
    $$ g_{1,1}(x,y,z) = x^2 - y^2 - x + y - 2z. $$
and since $g_{1,1}(-n,0,0) = n^2 + n$, at any time $t$ we have
    \begin{equation}
    \label{tasty}
    x^2 - y^2 = n^2 + n + x - y + 2z < n^2 + n.
    \end{equation}
For the second inequality, suppose first that $x+y \leq -n$.  Then
    $$ (x-1)^2 - (y-1)^2 = (x-y)(x+y-2) \geq -n(-n-2) > n^2. $$
On the other hand if $x+y > -n$, then $n+y+2z > 2z - x \geq x$, so from
(\ref{tasty}),
    $$ (x-1)^2 - (y-1)^2 = x^2 - y^2 - 2x + 2y
                 = n^2 + n - x + y + 2z > n^2. \qed $$
\renewcommand{\qedsymbol}{}
\end{proof}


\begin{prop}
Suppose that $a,b,n \in \N$ are such that $a^2 + n^2 = b^2$.  In the
situation of Proposition~\ref{hyperbolas}, for $t=a+b-n$, we have
$x(t)=-b$ and $y(t)=a$.
\end{prop}

\begin{proof}
By Lemma~\ref{dumb}, there exists some $t$ such that $y(t)=a$.  Now by the
first inequality of Proposition~\ref{hyperbolas}, for this value of $t$ we
have
    $$ x^2 < n^2 + n + a^2 = b^2 + n < b^2 + b < (b+1)^2, $$
so $x \geq -b$.  Now by the second inequality of
Proposition~\ref{hyperbolas}
    $$ (x-1)^2 > n^2 + (a-1)^2 = b^2 - 2a + 1 > (b-1)^2, $$
so $x \leq 1-b$.  Thus either $x=-b$ or $x=1-b$.  In the former case, the
proof is complete; in the latter case, we will show that $x(t+1)=-b$.
Indeed, if this were not the case, we would have $x(t+1)=1-b$ and
$y(t+1)= a+1$, but then
    $$ (x(t+1)-1)^2 - (y(t+1)-1)^2 = a^2 - b^2 = n^2, $$
\noindent contradicting the second inequality of
Proposition~\ref{hyperbolas}.
\end{proof}

\section{\large Higher-dimensional analogues}

\subsection{\normalsize The center of mass}

Recall that in dimension two and higher, the rotor-router model
depends on a cyclic ordering of the set of cardinal directions
$E_d = \{\pm e_1, \dots, \pm e_d\}  \subset \Z^d$.  Considering
$E_d$ as the set of vertices of a regular octahedron in $\R^d$, if
two orderings $\leq$ and $\leq'$ differ by an octahedral symmetry
of $E_d$, the resulting rotor-router models will differ by this
same symmetry.  By applying suitable reflections any ordering can
be transformed into one satisfying
 \begin{eqnarray}
        \label{conventions}
& & \text{\em (i) } ~ e_1 < e_2 < \dots < e_n;
~~~~~~~~~~~~~~~~~~~~~~~~~~~~~~~~~~ \nonumber\\
& & \text{\em (ii) } ~ e_i < -e_i, ~~ i=1,\dots,n.
~~~~~~~~~~~~~~~~~~~~~~~~~~~~~~~~~~
        \end{eqnarray}
For the remainder of this section, we fix an ordering $\leq$
satisfying (\ref{conventions})(i) and (ii).  Call a set $C \subset
[1,d]$ a {\it coclique} if the intervals $\{\epsilon \in E_d : e_i
< \epsilon \leq -e_i \}$, $i \in C$ are disjoint.

\begin{theorem}
\label{centerofmass}
The center of mass of the set of occupied sites is confined to the unit
cube $\{x \in \R^d : 0 \leq x_i \leq 1$\}.  Moreover, if $C \subset [1,d]$
is a coclique, then the center of mass lies below the hyperplane $\sum_{i
\in C} x_i = 1$.
\end{theorem}

\begin{proof}
For $\epsilon \in E_d$, let $T(\epsilon)$ be the total number of steps
taken in the direction $\epsilon$ by the first $n$ deposited particles.
Then $\frac{1}{n} \left[ T(e_i)-T(-e_i) \right]$ is the $i$-th coordinate
of the center of mass of the set of occupied sites after $n$ particles
have been allowed to equilibrate.  If the rotor $r(x)$ at the site $x \in
\Z^d$ satisfies $e_i < r(x) \leq -e_i$, then the site $x$ has ejected one
more particle in the direction $e_i$ than in the direction $-e_i$;
otherwise, $x$ has ejected equally many particles in the two directions.
Hence
    $$ 0 \leq T(e_i)-T(-e_i) \leq n, $$
and dividing by $n$, we conclude that the center of mass is confined to
the unit cube.

If $C \subset [1,d]$ is a coclique, then every occupied site $x$
satisfies at most one of the inequalities $e_i < r(x) \leq -e_i$, $i \in
C$, and so
    $$ 0 \leq \sum_{i \in C} \left[ T(e_i) - T(-e_i) \right] \leq n, $$
and dividing by $n$ gives the desired inequality.
\end{proof}

For example, suppose that $d=3$ and $\leq$ is the ordering $e_1 < e_2 <
-e_1 < -e_2 < e_3 < -e_3$.  Then the sets $\{1,3\}$ and $\{2,3\}$ are
cocliques, so the center of mass is confined to the portion of the unit
cube lying below the planes $x+z=1$ and $y+z=1$.

\subsection{\normalsize Progress toward circularity}

We conjecture that the limiting shape of the rotor-router model,
like that of IDLA, is a Euclidean ball in $\R^d$.  In this
section, we prove a much weaker, but analogous result,
Theorem~\ref{disc}.

\noindent The discrete Laplacian $\Delta F$ of a function $F : \Z^d
\rightarrow \R$ is given by
  $$ \Delta F(x) = \frac{1}{2d} \sum_{\epsilon \in E_d} F(x+\epsilon)
            -F(x). $$
If $\Delta F(x) = 0$, then $F$ is said to be harmonic at $x$.

Fixing an ordering of $E_d$, write $E_d = \{\epsilon_i\}_{i=1}^{2d}$ with
$\epsilon_1 < \dots < \epsilon_{2d}$.  Let $H_m(x)$ be the total number
of times the site $x \in \Z^d$ is been visited by the first $m$ deposited
particles.  The following result shows that $H_m$ is approximately
harmonic away from the origin.

\begin{lemma}
\label{harmonic}
If $x \neq {\bf 0}$, then $\Delta H_m(x)$ is bounded independent of $m$
and $x$.  Specifically,
    \begin{equation}
    \label{harmonicbounds}
    -d + \frac{3}{2} - \frac{1}{2d} \leq \Delta H_m
                    \leq d + \frac{1}{2}.
    \end{equation}
\end{lemma}

\begin{proof}
Every time a particle visits the site $x$, it comes from one of the neighboring
sites $x-\epsilon$, $\epsilon \in E_d$.  When the site $x-\epsilon_i$ is
first visited, the particle stays there, and thereafter, the $k$-th
particle to visit $x-\epsilon_i$ is routed to $x$ if and only if $k \equiv
i$ (mod $2d$).  The total number of routings from $x-\epsilon_i$ to $x$
after $n$ particles have been deposited is then at least $\frac{1}{2d}
\left[ H_m(x-\epsilon_i)-i \right]$ and, if $i<2d$, at most $a_i =
\frac{1}{2d} \left[ H_m(x-\epsilon_i)+2d-1-i \right]$.  In the case that
$i=2d$, we have $a_{2d} = \frac{1}{2d} \left[ H_m(x-\epsilon_{2d})-1
\right]$, whereas if the site $x-\epsilon_{2d}$ has not yet been visited,
then certainly no routings from $x-\epsilon_{2d}$ have taken place, so the
number of routings is actually $a_{2d} + \frac{1}{2d}$.  Summing the
contribution from each $x-\epsilon_i$, we obtain
    \begin{eqnarray*}
    & & -\frac{2d(2d+1)}{4d} + \frac{1}{2d}
        \sum_{\epsilon \in E_d} H_m(x-\epsilon) \leq H_m(x) \\
~&~&~~~~~~~~~~~~~~~~~~~~~~~~~~~
    \leq \frac{(2d-1)(2d-2)}{4d} + \frac{1}{2d}
        \sum_{\epsilon \in E_d} H_m(x-\epsilon).
    \end{eqnarray*}
and this reduces to (\ref{harmonicbounds}).
\end{proof}

There is a unique function $G$ on $\Z^d$ that is harmonic way from the
origin and satisfies $G(0)=0$ and $\Delta G(0)=-1$.  This $G$ is
called the discrete harmonic Green's function.  Unlike its continuous
counterpart, the discrete Green's function does not have a simple closed
form, and is given in two dimensions by an elliptic integral
\cite{Mangad}.  However, the discrete Green's function does have the same
asymptotics as its continous counterpart: in dimension $3$ and higher $G$
is asymptotic to a constant times $r^{2-d}$, and in dimension $2$ it is
asymptotic to a constant times $\log \frac{1}{r}$ (cf. \cite{Lawler96,
Mangad}).

A first attempt at a proof of the circularity conjecture might run as
follows.  To show that the rescaled set of occupied sites converges to a
ball, it would be enough to show that the function $H_m$ is sufficiently
radially symmetric, i.e. to bound $|H_m(x)-H_m(y)|$ in terms of $||x||^2
- ||y||^2$.  Due to its asymptotics, Green's function is itself
approximately radially symmetric, and given Lemma~\ref{harmonic}, we
might expect that the function $H_m$ should approximately coincide with
a suitable scaling and translation of $G$, namely $H_m({\bf 0}) - \Delta
H_m({\bf 0}) G$.  In two dimensions, however, we immediately encounter the
problem that $H_m$ is bounded below by zero, while $G \sim \frac{1}{2\pi}
\log \frac{1}{r}$ is not bounded below.  Indeed, the data suggest that $G$
is an excellent approximation to $H_m$ near $r=0$, but a bad approximation
when $H_m$ is close to zero; see figure~\ref{green'sfunction} and the
discussion in section~5.  Unfortunately, it is precisely when $H_m$ is
close to zero that we need a good approximation.

We might expect these difficulties to disappear in dimensions three and
higher, since $r^{2-d}$ is bounded below by zero.  Perhaps surprisingly,
however, Green's function does not appear to be a very good approximation
in the three dimensional case, either; see figure~\ref{3dgreen's} in
section~5.

To avoid the inaccuracies of the Green's function approximations, we will
take a somewhat different approach.  For the sake of simplicity, we treat
only the two-dimensional model.  Fix a map $r : \Z^2 \rightarrow E_2$
indicating the rotor direction at each point in the plane.  Imagine
now that several particles are simultaneously deposited at different points
in the plane; write $\nu(x)$  for the number of particles at the site $x$.
If $x$ is a point with $\nu(x) > 1$, let $x(r,\nu)$ denote the
configuration obtained by routing one particle from $x$ to the neighboring
site $x+r(x)$, and then changing the direction of the rotor $r(x)$ as
dictated by the ordering $\leq$.  A finite sequence of steps $(x) = (x_1,
x_2, \dots x_k)$ is said to be {\it terminating} if the configuration
$(r',\nu') = x(r,\nu) = (x_k x_{k-1} \dots x_1)(r,\nu)$ is such that
$\nu' \leq 1$ everywhere.  The following result is a special case
of Proposition 4.1 of \cite{DF}, but we include a proof here for the sake
of completeness.

\begin{prop}
\label{abelian}
Given a configuration of rotors $r_0:\Z^2 \rightarrow E_2$, and a map
$\nu_0 : \Z^2 \rightarrow \N$ indicating the number of particles at each
point in the plane, if $(x) = (x_1, \dots, x_k)$ and $(y) = (y_1, \dots,
y_l)$ are any two terminating sequences of steps, then $k=l$ and the
resulting configurations $x(r_0,\nu_0)$ and $y(r_0,\nu_0)$ are identical.
\end{prop}

\begin{proof}
It suffices to show that $(x)$ is a permutation of $(y)$.  If this were
not the case, then reversing the roles of $(x)$ and $(y)$ if necessary,
there exists $j$ such that the sequence $(x') = (x_1, \dots, x_{j-1})$ is
a permutation of a subsequence of $(y)$, but $(x'') = (x_1, \dots, x_j)$
is not.  Then $x_j$ occurs with the same multiplicity in $(x')$ and $(y)$,
while every $x_i \neq x_j$ occurs with at most the same multiplicity.  Now
setting $y(r_0,\nu_0) = (r, \nu)$ and $x'(r_0,v_0) = (r',\nu')$, it
follows that $\nu'(p) \leq \nu(p)$.  But $(y)$ is a terminating sequence,
so $\nu(p) \leq 1$, and hence $p$ is not a legal step after the sequence
$(x_1, \dots, x_j)$, a contradiction.
\end{proof}

Consider now the following procedure.  First, $m$ particles are deposited
simultaneously at the origin, and at each step thereafter, all but one of
the particles at each occupied site are routed to neighboring sites, until
there is at most one particle at each site.  By Proposition~\ref{abelian},
this procedure is guaranteed to terminate, and to give the same final
configuration $(r,\nu)$ our original model, in which the particles were
deposited one by one.  Notice that at each step in this new procedure,
each site ejects approximately equally many particles to each of its
neighbors; letting $H_{m,n}(x)$ be the number of particles at the site $x$
after $n$ steps, we see that $H_{m,n} \approx H_{m,n-1} + \Delta
H_{m,n-1}$, so our procedure has the effect of approximately iterating the
operator $\Delta + \text{Id}$.  Our next two lemmas convert this
observation into a precise estimate.

\begin{lemma}
\label{B_n}
Let $B_0(x,y) = \delta_{x0} \delta_{y0}$, where $\delta$ denotes the
Kronecker delta.  Let $B_n = B_{n-1} + \Delta B_{n-1}$, $n \geq 1$.  Then
    \begin{equation}
    \label{binomial}
    B_n(x,y) = 4^{-n} \binomial{n}{\frac12 (n+x+y)}
              \binomial{n}{\frac12 (n+x-y)},
    \end{equation}
where we adopt the convention that $\binomial{n}{k}=0$ if $k \notin \N$ or
$k>n$.
\end{lemma}

\begin{proof}
Induct on $n$.  Writing $u = \frac12 (n+x+y)$, $v = \frac12 (n+x-y)$, we
have by the inductive hypothesis
$$B_{n-1}(x-1,y) = 4^{1-n}\binom{n-1}{u-1} \binom{n-1}{v-1};$$
$$B_{n-1}(x,y-1)= 4^{1-n}\binom{n-1}{u-1} \binom{n-1}{v};$$
$$B_{n-1}(x,y+1)= 4^{1-n}\binom{n-1}{u} \binom{n-1}{v-1};$$
$$B_{n-1}(x+1,y)= 4^{1-n}\binom{n-1}{u} \binom{n-1}{v};$$
whence
    \begin{eqnarray*}
    B_n(x,y) &=& 4^{-n} \left[ \binomial{n-1}{u-1}
+ \binomial{n-1}{u} \right] \binomial{n-1}{v-1} \\ &~&~~~~~~~~~~~~~~~~~~~~
           + 4^{-n} \left[ \binomial{n-1}{u-1}
+ \binomial{n-1}{u} \right] \binomial{n-1}{v} \\
         &=& 4^{-n} \binomial{n}{u} \binomial{n}{v}. \qed
    \end{eqnarray*}
\renewcommand{\qedsymbol}{}
\end{proof}

\begin{lemma}
\label{3n}
$|H_{m,n} - mB_n| \leq 3n$.
\end{lemma}

\begin{proof}
By definition, $F_{m,0} = m \delta_{x0} \delta_{y0} = m B_0$.  We now
show by induction on~$n$ that
    \begin{equation}
    \label{inductivestep}
    -3 \leq H_{m,n} - (\Delta + \text{Id}) [H_{m,n-1}] \leq 3.
    \end{equation}
On the $n$-th step, depending on the direction of the rotor
$r(z+\epsilon)$ and the congruence class of $H_{m,n-1}(z+\epsilon)$ (mod
$4$), a point $z \in \Z^2$ receives between $\frac14 \left[
H_{m,n-1}(z+\epsilon) - 4 \right]$ and $\frac14 \left[
H_{m,n-1}(z+\epsilon) + 2 \right]$ particles from each neighboring site
$z + \epsilon$.  Moreover, the site $z$ ejects all but one of its own
particles, leaving
    $$-4 + (\Delta + \text{Id}) [H_{m,n-1}](z) + 1 \leq H_{m,n}(z)
        \leq 2 + (\Delta + \text{Id}) [H_{m,n-1}](z) + 1, $$
and (\ref{inductivestep}) follows.
\end{proof}

We're now ready to prove the promised weak circularity result.  We will
show that after $m$ particles are deposited and allowed to equilibrate,
every site in a disc centered at the origin of radius proportional to
$m^{1/4}$ is occupied.  We will make use of Stirling's inequality
    \begin{equation}
    \label{stirling's}
    \binomial{2n}{n} > \frac{4^n}{e^2\sqrt{n}}.
    \end{equation}
and the fact that for any $\epsilon>0$ there exists $N(\epsilon)$
such that for all $t > N(\epsilon)$,
    \begin{equation}
    \label{e^a}
    \left( 1-\frac{a}{t} \right)^t  > e^{-a-\epsilon}
    \end{equation}

\begin{theorem}
\label{disc}
Let $\epsilon>0$.  In either the type 1 or type 2 rotor-router model
on $\Z^2$, if $m$ is taken sufficiently large, then after $m$
particles are deposited and allowed to equilibrate, every site in the
open disc of radius $r_0 = \left( \frac{8m}{3e^{6+\epsilon}}
\right)^{1/4}$ centered at the origin is occupied.
\end{theorem}

\begin{proof}
Let $(x,y) \in \Z^2$ and suppose first that $x \equiv y \equiv 0$
(mod $2$).  Write $u = \frac12 |x+y|$, $v = \frac12 |x-y|$. Let $r
= \sqrt{x^2 + y^2} = \sqrt{2u^2+2v^2}$ and let $n = \frac12 r^2$.
Note that $n$ is an even integer.  By Lemmas~\ref{B_n}
and~\ref{3n} and Stirling's inequality (\ref{stirling's}), for any
$m \in \N$ we have
    \begin{eqnarray}
    H_{m,n}(x,y)  &\geq&  m 4^{-n} \binomial{n}{u + n/2}
        \binomial{n}{v + n/2} - 3n \nonumber \\
&=& m 4^{-r^2/2} \binomial{r^2/2}{r^2/4}^2 \cdot
    \prod_{i=1}^{u} \frac{\frac{r^2}{4}-i}{\frac{r^2}{4}+i}
    \prod_{j=1}^{v} \frac{\frac{r^2}{4}-j}{\frac{r^2}{4}+j}
    - \frac32 r^2 \nonumber \\
&>& \frac{4m}{e^4 r^2}
    \prod_{i=1}^{u}
        \frac{\frac{r^2}{4}-u-1+i}{\frac{r^2}{4}+i}
    \prod_{j=1}^{v}
        \frac{\frac{r^2}{4}-v-1+j}{\frac{r^2}{4}+j}
    - \frac32 r^2 \nonumber \\
&=& \frac{4m}{e^4 r^2}
    \prod_{i=1}^{u}
        \left(1-\frac{4(u+1)}{r^2+4i} \right)
    \prod_{j=1}^{v}
        \left(1-\frac{4(v+1)}{r^2+4j} \right)
    - \frac32 r^2 \nonumber \\
&>& \frac{4m}{e^4 r^2}
    \left( 1 - \frac{4(u+1)}{r^2} \right)^{u}
    \left( 1 - \frac{4(v+1)}{r^2} \right)^{v}
    - \frac32 r^2
    \label{e!}
    \end{eqnarray}
We now use the inequality (\ref{e^a}) in the cases
    \begin{eqnarray*}
& & a_1 = \frac{4u(u+1)}{r^2}, ~~~t_1 = u;
~~~~~~~~~~~~~~~~~~~~~~~~~~~~~~~~~~ \\
& & a_2 = \frac{4v(v+1)}{r^2}, ~~~t_2 = v.
~~~~~~~~~~~~~~~~~~~~~~~~~~~~~~~~~~
    \end{eqnarray*}
Let $\delta$ be such that $2 \delta + \frac{2}{N(\delta)} =
\epsilon$. For $u,v > N(\delta)$ and $r < \left( \frac{8m}{3
e^{6+\epsilon}} \right)^{1/4}$, we have from (\ref{e!}),
    \begin{eqnarray}
    \label{almost}
H_{m,n}(x,y) &>& \frac{4m}{e^4 r^2}
    \exp \left( -\frac{4u^2 + 4v^2 + 4u + 4v}{r^2}
        - 2 \delta \right)
    - \frac32 r^2 \nonumber \\
&=& \frac{4m}{e^{6+2\delta} r^2}
    \exp \left( -\frac{2u+2v}{u^2+v^2} \right)
    - \frac32 r^2 \nonumber \\
&>& \frac{4m}{e^{6+2\delta} r^2}
    \exp \left( -\frac{2\sqrt{2}}{\sqrt{u^2+v^2}} \right)
    - \frac32 r^2 \nonumber \\
&>& \frac{4m}{e^{6+2\delta} r^2}
    \exp \left( -\frac{2}{N(\delta)} \right)
    - \frac32 r^2 \\
&=& \frac{3}{2r^2}
    \left[ \frac{8m}{3e^{6+\epsilon}} - r^4 \right]
> 0. \nonumber
    \end{eqnarray}
It remains to consider the case that $u$ or $v$ is $\leq
N(\delta)$. First, let $\beta = \left( \frac{4}{1-e^{-\delta}}
\right)^{1/2}$, and take $m$ large enough so that $m B_n(x,y) >
3n$ for all pairs $x,y$ for which $u,v \leq \beta N(\delta)$.  By
Lemma~\ref{3n}, $H_{m,n}(x,y) > 0$ for all such pairs.  Now
suppose one of $u,v$ is $> \beta N(\delta)$. Taking $y$ to $-y$ if
necessary, we can assume $u \leq N(\delta)$, $v > \beta
N(\delta)$.  In particular, this implies
    $$ \frac{4u(u+1)}{r^2} < \frac{4u^2}{v^2}
        < \frac{4}{\beta^2} = 1 - e^{-\delta}, $$
so from (\ref{e!}) we obtain
    \begin{eqnarray*}
H_{m,n}(x,y) &>& \frac{4m}{e^4 r^2}
    \left(1 - \frac{4u(u+1)}{r^2} \right)
    \exp \left( -\frac{4v(v+1)}{r^2} - \delta \right)
    - \frac32 r^2 \\
&>& \frac{4m}{r^2} \exp \left(-4 - 2\delta -
        \frac{2v^2}{u^2 + v^2} - \frac{2v}{u^2 + v^2} \right)
    - \frac32 r^2 \\
&>& \frac{4m}{r^2} \exp \left(-6 - 2\delta - \frac{2}{N(\delta)} \right)
    - \frac32 r^2
    \end{eqnarray*}
and we recover equation (\ref{almost}).

In the cases $x \equiv y \equiv 1$ (mod $2$) and $x \not\equiv y$
(mod $2$) the proof is similar, taking $n = \frac{r^2-2}{2}$ and
$\frac{r^2-3}{2}$, respectively.
\end{proof}

\section{\large Conjectures}

\subsection{\normalsize The Sturmian region}

\begin{figure}
\centering
\includegraphics[scale=.9]{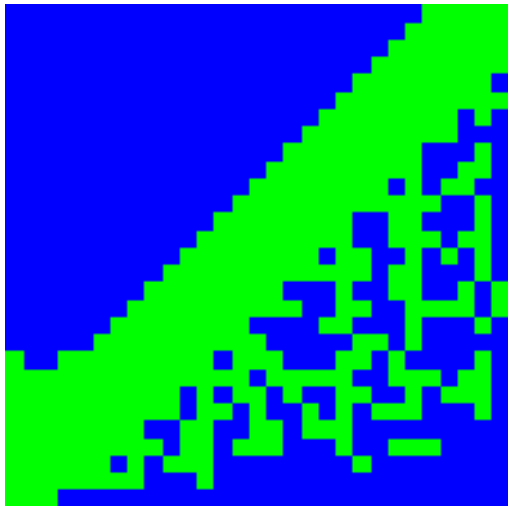}
\caption{\footnotesize Plot of the set of pairs $(r,s)$, $1 \leq
r,s < 30$, for which $w_{r,s)}$ is Sturmian out to $10^7$ places;
$r$ increases from left to right and $s$ increases from bottom to
top.  The square of area $1$ centered at $(r,s)$ is shaded green
(light gray) if $w_{r,s}$ is Sturmian, and blue (dark gray) if
$w_{r,s}$ is not Sturmian.} \label{sturmianregion}
\end{figure}

Figure \ref{sturmianregion} shows the set of pairs $(r,s)$ with $1
\leq r,s < 30$ for which the binary word $w_{r,s}$ is Sturmian in
the first ten million terms. In each of these cases, $w_{r,s}$
takes the form (\ref{floor}) with
    \begin{equation}
    \alpha = \frac{\sqrt{s}}{\sqrt{r}+\sqrt{s}}, ~~~~~
    \beta = \frac{\alpha-1}{r}+\frac12.
    \label{alphabeta}
    \end{equation}
We make the obvious conjectures.

\begin{conjecture} ~ \\
\label{obvious}
\indent \indent (i) For $(r,s)$ in the diagonal stripe $-4 \leq r-s
\leq 3$, with the single exception of the case $r=4$, $s=1$ the
sequence $w_{r,s}$ is Sturmian. \\
\indent \indent ~(ii) If $(r,s)$ is just below the stripe, i.e.
$r-s = 4$, then $w_{r,s}$ is Sturmian if and only if $r$ is even. \\
\indent \indent ~(iii) If $r-s < -4$, then with finitely many
exceptions, $w_{r,s}$ is not Sturmian.
\end{conjecture}

It is likely that (i) and the positive direction of (ii) can be
proved in the same way as Proposition~\ref{sturmian}; the
computations become quite extensive, however.  Proving the
negative direction of (ii), as well as (iii), may be trickier.  It
is not hard to show that if $w_{r,s}$ is Sturmian, then its slope
$\alpha$ must be given by (\ref{alphabeta}); indeed, from Theorem
\ref{boundedness} we have $\frac{|x|}{y} \rightarrow
\sqrt{\frac{r}{s}}$ as $t \rightarrow \infty$, and so the
proportion of $u \leq t$ for which $f(\sigma(u)) = f^-(\sigma(u))$
approaches
    $$\frac{\frac{|x|}{r}}{\frac{|x|}{r}+\frac{y}{s}}
    = \frac{\frac{|x|}{ry}}{\frac{|x|}{ry}+\frac{1}{s}} \longrightarrow
      \frac{\frac{1}{\sqrt{rs}}}{\frac{1}{\sqrt{rs}}+\frac{1}{s}}
    = \frac{\sqrt{s}}{\sqrt{r}+\sqrt{s}}, $$
and we recover (\ref{alphabeta}).  Likewise, from (\ref{bounds})
one can deduce that if $w_{r,s}$ is Sturmian, its intercept
$\beta$ must be given by (\ref{alphabeta}).

The fact that $w_{r,s}$ appears to be Sturmian for so many pairs
$(r,s)$ suggests that an exact description of the triples
$f_{r,s}^n(0,0,0)$ for general~$r$ and~$s$ may be within reach.
The subword complexity (number of factors of length~$n$) of
$w_{r,s}$ would be worth investigating in this connection, since
Sturmian words are of minimal complexity among aperiodic words
(cf. \cite{BS}).

\subsection{\normalsize Green's function estimates}

\begin{figure}
\centering
\includegraphics[scale=.33, angle=270]{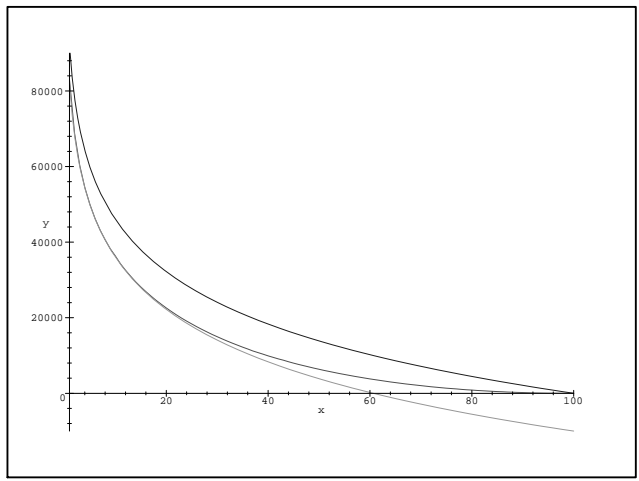}
\caption{\footnotesize The function $H_m(r,0)$ (middle curve) for
the type 1 rotor-router model on $\Z^2$, plotted against the
functions $F(r)$ (lower curve) and $\tilde{F}(r)$ (upper curve) of
equation (\ref{F}). Here $m = \floor{10000\pi}$, so that $H_m$ has
its root near $r=100$.} \label{green'sfunction}
\end{figure}

The discrete Green's function in two dimensions has the
asymptotics (cf. \cite{Mangad})
    $$ G(r) = \frac{1}{2\pi} \left( \log r + \frac32 \log 2
        + \gamma + O(\frac{1}{r^2}) \right). $$
Set $\tilde{G}(r) = \frac{1}{2\pi} \left (\log r + \frac32 \log 2
+ \gamma \right)$.  As in section 4.2, let $H_m(x,y)$ denote the
total number of times the point $(x,y) \in \Z^2$ is visited by the
first $m$ deposited particles.  As discussed in section 4.2, we
might expect $H_m$ to coincide closely with the function $H_m(0,0)
- \Delta H_m(0,0)G'$.  Observe that every time a particle visits
the origin, it has either been newly placed there, or it has come
from one of the four adjacent sites $(\pm 1, 0)$, $(0, \pm 1)$; so
by the same argument used in the proof of Lemma~\ref{harmonic}, we
conclude that $\Delta H_m(0,0)$ is approximately~$-m$.
Figure~\ref{green'sfunction} plots $H_m(r,0)$ against the two
functions
    \begin{equation}
    \label{F}
    F = H_m(0,0) - m\tilde{G}, ~~~~~
       \tilde{F} = H_m(0,0) - m\left(\tilde{G} - \frac12 \right).
    \end{equation}
Since we constructed $F$ to coincide with $H_m$ near zero, and $F$
is unbounded below while $H_m \geq 0$, it is not too surprising
that $F$ gives a good approximation to $H_m$ near zero, but a bad
approximation near the root of $H_m$.  Interestingly, however, the
root of $\tilde{F}$ seems to coincide very closely with that of
$H_m$.

\begin{figure}
\centering
\includegraphics[scale=.35, angle=270]{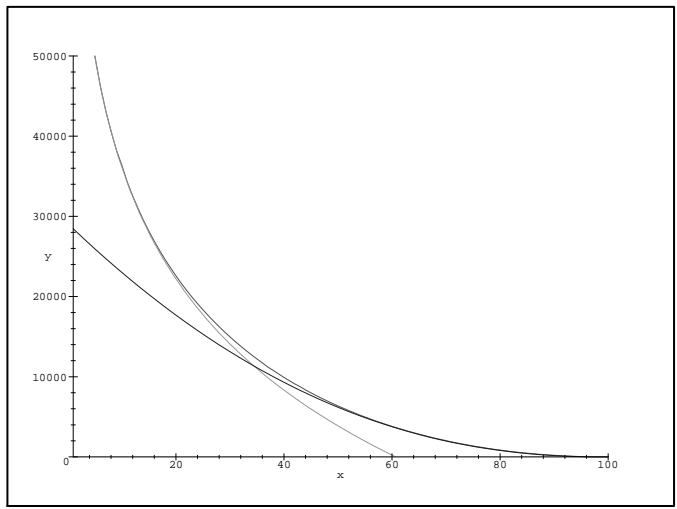}
\caption{\footnotesize Near its root $r_0$, $H_m(r,0)$ falls off like
$|r-r_0|^\lambda$, where $\lambda$ varies only very slightly with $m$.
Shown here is $H_m$ (upper curve) with $m=\floor{10000\pi}$,
alongside the curves $H(r)$ and $|r-r_0|^\lambda$ for $\lambda=2.232$.}
\label{boundatroot}
\end{figure}

\begin{figure}
\centering
\includegraphics[scale=.33, angle=270]{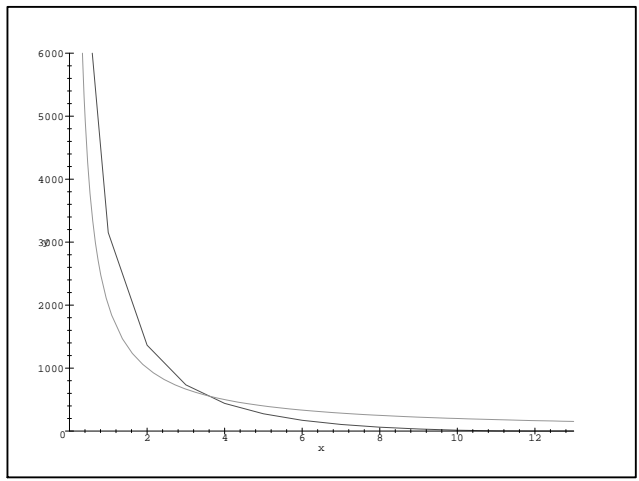}
\caption{\footnotesize The function $H_m(r,0,0)$ for a three-dimensional
rotor-router model (darker curve), plotted against $m/r$, the dominant
term of Green's function on $\Z^3$.  Here $m=\floor{(4/3)\pi\cdot 15^3}$,
so that $H_m$ has its root near $r=15$.}
\label{3dgreen's}
\end{figure}

Figure~\ref{boundatroot} shows $H_m$ plotted against $H$ and the
function $|r-r_0|^\lambda$, where $r_0 \approx 100$ is the root of
$H_m$ and $\lambda = 2.232$; this latter function gives a very
good approximation to $H_m$ near its root.

Since Green's function is bounded below on $\Z^n$ for $n \geq 3$,
we might expect it to give better approximations to $H_m$ in
higher dimensions than it does in dimension two.  For the most
part, however, these expectations do not seem to be bourne out.
Figure~\ref{3dgreen's} shows $H_m(r,0,0)$ for the
three-dimensional rotor-router model obtained from the ordering
$e_1 < -e_1 < e_2 < -e_2 < e_3 < -e_3$, plotted against $m/r$, the
principal term of the discrete Green's function in three
dimensions.

~ \\ {\bf Acknowledgments.} The author would like to thank Jim Propp for 
his invaluable advising, and for providing the initial motivation for 
much of this work.  In addition, Adam Kampff helped write the C code 
for generating many of the figures reproduced here.  Henry Cohn and Anna 
Salamon provided helpful comments on earlier drafts of this paper.

\end{document}